\documentclass[12pt]{amsart}      
\usepackage{amsfonts}
\usepackage{mathrsfs}
\usepackage{amssymb}
\newtheorem{proposition}{Proposition}[section]

\theoremstyle{definition}

\theoremstyle{remark}

\numberwithin{equation}{section}

\begin{document}
\baselineskip=20pt

\title{Principal Realization of the Yangian $Y(\frak g\frak l (n))$}

\author{Cheng-Ming Bai}
\address{Bai: Chern Institute of
Mathematics and LPMC, Nankai University, Tanjin 300071, China}
\email{baicm@nankai.edu.cn}

\author{Mo-Lin Ge}
\address{Ge: Chern Institute of
Mathematics, Nankai University, Tanjin 300071, China}
\email{geml@nankai.edu.cn}

\author{Naihuan Jing$^*$}
\address{Jing: Department of Mathematics, North Carolina State University,
Raleigh, NC 27695-8205, U.S.A.}
\email{jing@math.ncsu.edu}

\def\shorttitle{Principal Realization of the Yangian $Y(\frak g\frak l (n))$}

\begin{abstract} Motivated to simplify the structure of
tensor representations we give a new set of generators for the
Yangian $Y(sl(n))$ using the principal realization in simple Lie
algebras. The isomorphism between our new basis and the standard
Cartan-Weyl basis is also given. We show by example that the
principal basis simplifies the Yangian action significantly in the
tensor product of the fundamental representation and its dual.

\end{abstract}


 \thanks{$^*$ Corresponding
author.}

\maketitle

\section{Introduction}

Loop algebras $L(\frak g)$ and affine Lie algebras $\hat {\frak g}$
are important generalizations of the finite dimensional simple Lie
algebra $\frak g$. Their structures and representations are usually
studied with the help of particular presentations of the Lie
algebra. There are two important realizations for the Lie algebra
$\frak g$ as well as $L(\frak g)$ and $\hat {\frak g}$: the
Cartan-Weyl presentation and the principal realization \cite{K}. For
example, in the case of $\frak g=\frak g\frak l(n)$ the Cartan-Weyl
basis consists of the usual unit matrices $\{ E_{ij}\}_{1\leq
i,j\leq n}$ while the principal basis elements $\{ A_{ij}\}_{i,j\in
{\mathbb Z}/{n{\mathbb Z}}}$ are Toeplitz matrices, which are
diagonal-constant or those in which each descending diagonal from
left to right is constant. The principal realization played a
pivotal role in constructing the first example of the vertex
representation of the affine Lie algebra and opened a new chapter
for the conformal field theory.

Yangians are certain quantum deformations of the loop algebras
associated to the simple Lie algebras. In the study of the rational
solution to the Yang-Baxter equation Drinfeld defined these Hopf
algebras as special symmetry algebras \cite{C, CP}, and they have
been actively studied from various contexts (see the recent
comprehensive monograph \cite{M} and also \cite{AMR, BR, BK, I}). In
\cite{D1} V.G. Drinfeld introduced the Yangian for an arbitrary
simple Lie algebra as an associative algebra generated by $\{
I_\alpha, J_\alpha\}$ which are analog of the loop algebra and then
found new realization for the Yangian in \cite{D2}. The realization
of $Y(\frak g\frak l (n))$ via the RTT-equation goes back to
Faddeev's school \cite{KS}. The generators $\{T_{ij}^{(n)}\}$ in the
FRT formulation is usually better for studying properties resembling
the matrix operations such as the quantum determinant \cite{M},
while Drinfeld's new basis plays a crucial role in studying finite
dimensional representations (cf. \cite{CP} and \cite{M}). However
the second part of Drinfeld generators called $J_\alpha$ operators
are often difficult to compute in various representations due to a
symmetric sum in the construction. In \cite{BGX} an example is found
that in certain change of basis the action of the $J_{\alpha}$
operators seem to be simplified. It is interesting that this example
was found in the application of Yang-Baxter $R$-matrix to maximally
entangled states in the context of the newly active field of quantum
computation. This raises the question whether there exists another
basis fitted for this purpose.

In this paper we
 shall introduce
another basis or presentation for the Yangian $Y(\frak g\frak l(n))$
and $Y(\frak s\frak l(n))$ using the principal gradation. Our
construction is motivated by two reasons. The first is to establish
a mathematical model to answer our question in the last paragraph
and seek a better way to study representations of the Yangian. There
seems to be some interesting properties shared by our new generators
and in particular we exploit the additive structure in the index set
$I=\{ 0, 1,\cdots, n-1\}$ by using the discrete Fourier transform on
the finite abelian group $\mathbb Z_n$. It turns out that in the
principal basis some of the actions of the $J$ operators are
dramatically simplified, and we offer one example in the last
section to support this claim. The second motivation is the
connection with quantum computation. We notice that similar
construction of the principal basis has been used in several fast
algorithms in quantum computation through the discrete Fourier
transform. In this aspect we will show by an example that our
principal realization indeed has a close relationship with quantum
entanglement and it simplifies some of complicated computations
experienced in the standard Cartan-Weyl basis for the Yangian.

The paper is organized as follows. Section two reviews some of the
basic materials about principal basis for the general linear Lie
algebra. Section three discusses how to change the Cartan-Weyl basis
to another basis and introduces the principal realization for the
Yangian $Y(\frak g\frak l(n))$. We also study basic properties for
the principal basis and give the isomorphism between the Cartan-Weyl
basis and our principal basis. Section four studies the
representation of the Yangian in terms of the principal basis.
Finally in section five we apply our basis to the tensor product of
the fundamental representations of the Yangian $Y(sl(3))$ and show
that the principal basis can simplify the action of the Yangian in
this case.

\section{The general linear Lie algebra $\frak g\frak l(n)$}

Let $\frak g=\frak g\frak l(n)$ be the Lie algebra of $n\times
n$-complex matrices. The standard unit matrices $\{ E_{ij}\}$ form
the so-called Cartan-Weyl basis for $\frak g$ and
$$[E_{ik}, E_{lj}]=\delta_{kl}E_{ij}-\delta_{ij}E_{lk},\;\;\;i,j,k,l\in I.\eqno (2.1)$$
The enveloping algebra is given by the relation:
$E_{ik}E_{lj}=\delta_{kl}E_{ij}$. If we define ${\rm
deg}(E_{ij})=j-i$, then $\frak g$ becomes a graded Lie algebra
$$\frak g=\bigoplus_{i\in {\mathbb Z}/{n{\mathbb Z}}}{\frak g}_i,\;\;
[{\frak g}_i, {\frak g}_j]\subset {\frak g}_{i+j},\eqno (2.2)$$
where ${\frak g}_i=\{ x\in \frak g|{\rm deg} (x)=i\}$. For a matrix
$x=(x_{ij})$ we define the principal decomposition:
$$x=x_0+x_1+x_2+\cdots+x_{n-1},\;\;{\rm deg}(x_i)=i,\eqno (2.3)$$
where $x_i=(x_{k.k+i})_{k,l}$. For instance
$$x_1=\left(\begin{matrix} 0 & x_{12} &0 & \cdots & 0\cr
0 &0 & x_{23} & \cdots & 0\cr \vdots& \vdots &\vdots &\vdots
&\vdots\cr 0 &0 &0 &\cdots & x_{n-1,n}\cr x_{n1}& 0 &0 &\cdots
&0\cr\end{matrix}\right).$$

Consider $\sigma\in {\rm Aut}(\frak g)$ defined by
$$\sigma(E_{j,j+1})=\omega
E_{j,j+1},\;\;\sigma(E_{j+1,j})=\omega^{-1}E_{j+1,j},\;\;\omega=e^{\frac{2\pi
i}n},\eqno (2.4)$$ and extended to $\frak g$ linearly. Then ${\frak
g}_i=\{ x\in \frak g|\sigma(x)=\omega^i x\}$. Let
$$E=\sum_{i=0}^{n-1} E_{i,i+1}$$ then the
centralizer $Z(E)=\bigoplus\limits_{i=0}^{n-1}{\mathbb C} E^i$ is a
Cartan subalgebra of $\frak g\frak l(n)$ called the principal Cartan
subalgebra. With respect to the principal Cartan subalgebra, the
(principal) root spaces are defined by ${\frak g}_\beta=\{ x\in
\frak g|[h,x]=\beta(h)x,\;\;\forall\; h\in Z(E)\}$. In this case the
root vectors are given as follows.

For $i\in{\mathbb
Z}/{n{\mathbb Z}}$ let
$$A_i=(\omega^{ki})_{kl}=\left( \begin{matrix} \omega^i & \omega^i
&\cdots & \omega^i\cr \omega^{2i} & \omega^{2i}& \cdots
&\omega^{2i}\cr \vdots &\vdots&\vdots&\vdots\cr \omega^{ni}
&\omega^{ni} &\cdots &\omega^{ni}\cr\end{matrix}\right)$$ and we
denote its principal components by $A_{ij}$:
$$A_{i}=A_{i0}+A_{i1}+\cdots+A_{i,n-1},\eqno (2.5)$$
i.e. $A_{ij}=(\omega^{ki}\delta_{l-k,j})_{k,l}=\sum_k\omega^{ik}E_{k, j+k}$
 for $j\in {\mathbb
Z}/{n{\mathbb Z}}$ and $i=1,\cdots, n-1$.
 For completeness we include the Lie algebra structure under the
principal decomposition as follows.
$$
[E^k, A_i]=(\omega^{ki}-1)A_i, \quad [E^k,
A_{ij}]=(\omega^{ki}-1)A_{i, j+k}.\eqno(2.6)$$
$$[A_{ij},
A_{i'j'}]=(\omega^{ji'}-\omega^{j'i})A_{i+i',j+j'}, \eqno(2.7)
$$
which follows from the algebra structure:
$A_{ij}A_{i'j'}=\omega^{ji'}A_{i+i',j+j'}$.

Under the standard invariant form $(x|y)={\rm tr}(xy)$, the basis
elements are isotropic except $n|2i$ and $n|2j$:
$(A_{ij}|A_{i'j'})=n\omega^{-ij}\delta_{i,-i'}\delta_{j,-j'}$.

\section{The Yangian $Y(\frak g\frak l(n))$}

The Yangian $Y(\frak g \frak l(n))$ is defined by Drinfeld \cite{D1,
KS} as the complex associative unital algebra generated by
generators $T_{ij}^{(m)}$, $i,j\in\{ 1,\cdots, n\}, m=1,2,\cdots $
subject to the defining relations (cf. \cite{M}):
$$[T_{ij}^{(l+1)},T_{i'j'}^{(m)}]-[T_{ij}^{(l)},T_{i'j'}^{(m+1)}]
=T_{i'j}^{(l)}T_{ij'}^{(m)}-T_{i'j}^{(m)}T_{ij'}^{(l)},\eqno (3.1)$$
where $l,m=0,1,2\cdots$ and $T_{ij}^{(0)}=\delta_{ij}\cdot 1$. For
this reason the generators $T_{ij}^{(m)}$ are Cartan-Weyl type
generators. It is well known that $\{ T_{ij}^{(1)}\}\cup\{
T_{ij}^{(2)}\}$ are enough to span the whole Yangian $Y(\frak g\frak
l(n))$. The following result can be easily proved using dual bases.

{\bf Lemma 3.1}\quad Let $\{ e_i\}$ and $\{ e^i\}$ be a pair of dual
bases of $\frak g$. Then the $r$-matrix can be expressed as follows.
$$r=\sum e_i\otimes e^i.\eqno (3.2)$$
Moreover this expression is independent of the choice of the dual
bases.

Using the principal basis (2.7) it is easy to get the following.

{\bf Corollary 3.2}\quad The permutation matrix can be written as
$$P=\sum_{ij}E_{ij}\otimes E_{ji}=\sum_{k,l}\frac{\omega^{kl}}{n}
A_{kl}\otimes A_{-k,-l}.$$

Let $u$ be a formal variable and we consider the operators or
matrices in ${\rm End}(V\otimes V)[[u^{-1}]]$ as well as $Y(\frak
g\frak l (n))[[u^{-1}]]$. Let
$$R(u)=I-\frac{P}{u}\in {\rm End}(V\otimes V)[[u^{-1}]].$$

 The Yang-Baxter equation is an operator equation on $End(V^{\otimes 3})[[u^{-1}]]$
written as
$$R_{12}(u)R_{13}(u+v)R_{23}(v)=R_{23}(v)R_{13}(u+v)R_{12}(u),$$
where the indices specify the action of $R$ on respective components
of $V^{\otimes 3}$. The so-called $T$-matrix $T(u)=(T_{ij}(u))$ of
the Yangian $Y(\frak g\frak l(n))$ can incorporate the defining
relations (3.1) into a mater matrix identity as follows. Let
$$T_{ij}(u)=\delta_{ij} +\sum_{k=1}^\infty T_{ij}^{(k)}u^{-k}\in
Y(\frak g\frak l(n))[[u^{-1}]].$$ Then
$$T(u)=\sum_{i.j}T_{ij}(u)\otimes E_{ij} \in Y(\frak
g\frak l(n))\otimes {\rm End}(V). \eqno (3.3)$$

{\bf Proposition 3.3}\quad (\cite{D1, D2, KS})\quad The defining
relations of the Yangian can be written as
$$R(u-v)T_1(u)T_2(v)=T_2(v)T_1(u)R(u-v) \eqno (3.4)$$
and the coproduct is given by
$$\Delta(T(u))=T(u)\otimes T(u).$$

We now introduce a new set of generators, which are generalization
of the principal generators in $\frak g\frak l(n)$. For $k,l\in
{\mathbb Z}/{n{\mathbb Z}}$ we define elements
$$S_{kl}(u)=\sum_{m=0}^\infty S_{kl}^{(m)} u^{-m},\eqno (3.5) $$
and $S_{kl}^{(0)}=\delta_{k0}\delta_{l0} \cdot 1$.

We write
$$T(u)=\sum_{k,l\in
{\mathbb Z}/{n{\mathbb Z}}} S_{kl}(u)\otimes A_{kl}\in Y(\frak
g\frak l(n))\otimes {\rm End}(V).$$

The following result is a standard calculation in FRT formulation.

 {\bf Theorem 3.4}\quad The principal generators
$S_{kl}^{(m)}$ of the Yangian satisfy the following relations:
$$[S_{ij}(u),S_{kl}(v)]=\frac{1}{u-v}(\sum_{a,b}
\frac{\omega^{ib-bk-ab}}{n}S_{k+a,l+b}(u)S_{i-a,j-b}(v)
$$$$-\sum_{a,b}\frac{\omega^{ja-al+ab}}{n}S_{k+a,l+b}(v)S_{i-a,j-b}(u)).$$

{\bf Proof}. \quad From the FRT equation it follows that
$$(I-\frac{P}{u-v})(\sum_{i,j,k,l}S_{ij}(u)S_{kl}(v)A_{ij}\otimes
A_{kl})=(\sum_{i,j,k,l}S_{ij}(v)S_{kl}(u)A_{ij}\otimes
A_{kl})(I-\frac{P}{u-v}).$$ Using the principal decomposition of the
operator $P$ we have
\begin{eqnarray*}
[S_{ij}(u),S_{kl}(v)]&=& \frac{1}{u-v}(P\sum_{i,j,k,l}(A_{ij}\otimes
A_{kl}) S_{ij}(u)S_{kl}(v)-\sum_{i,j,k,l}(A_{ij}\otimes A_{kl})
S_{kl}(v)S_{ij}(u) \cdot P)\\
&=&\frac{1}{u-v}(\sum_{i,j,k,l,a,b}\frac{\omega^{ab+bi-bk}}{n}(A_{i+a,j+b}\otimes
A_{k-a,l-b})S_{ij}(u)S_{kl}(v)\\
&\mbox{}&
-\sum_{i,j,k,l,a,b}\frac{\omega^{ab+ja-la}}{n}(A_{i+a,j+b}\otimes
A_{k-a,l-b})S_{kl}(v)S_{ij}(v))\\
&=&
\frac{1}{u-v}\sum_{i,j,k,l,a,b}\frac{\omega^{-ab+bi-bk}}{n}(A_{ij}\otimes
A_{kl})S_{k+a,l+b}(u) S_{i-a,j-b}(v)\\
&\mbox{}&-\frac{1}{u-v}\sum_{i,j,k,l,a,b}\frac{\omega^{ab+ja-la}}{n}(A_{ij}\otimes
A_{kl})S_{k+a,l+b}(v) S_{i-a,j-b}(u).
\end{eqnarray*}
\hfill $\Box$

The defining relations of our principal realization can be written
componentwise as follows.
\begin{eqnarray*}
[S_{ij}^{(l+1)},S_{i'j'}^{(m)}]-[S_{ij}^{(l)},S_{i'j'}^{(m+1)}]&=&
\sum_{a,b\in{\mathbb Z}/{n{\mathbb Z}}}
\frac{\omega^{(i-i')b-ab}}{n}S_{i'+a,j'+b}^{(l)}S_{i-a,j-b}^{(m)}\\
&\mbox{}& -\sum_{a,b\in{\mathbb Z}/{n{\mathbb Z}}}
\frac{\omega^{(j-j')b+ab}}{n}S_{i'+a,j'+b}^{(m)}S_{i-a,j-b}^{(l)}
\end{eqnarray*}
where all lower indices are inside ${\mathbb Z}/{n{\mathbb Z}}$, and
$l,m\in {\mathbb N}=\{0,1,\cdots\}$.

{\bf Theorem 3.5.} The mapping $\displaystyle T_{ij}(u)\mapsto
\sum_k S_{k, j-i}(u)\omega^{ik}$ defines an isomorphism of two
presentations of the Yangian $Y({\mathfrak gl}(n))$. The inverse
mapping is given by
$$\displaystyle S_{kl}(u)\mapsto \sum_i
\frac{\omega^{-ki}}nT_{i+l,i}(u). \eqno(3.6)$$

{\bf Proof}. This can be simply shown by the master equations
(3.3-3.4). We prove the inverse isomorphism. Note that $\{ A_{kl}\}$
and $\{\frac{\omega^{kl}}nA_{-k, -l}\}$ are dual bases.
\begin{align*}
S_{kl}(u)&=(T(u)|A_{-k,-l})\frac{\omega^{kl}}n=
\sum_{ij}T_{ij}(u)(E_{ij}|A_{-k,-l})\frac{\omega^{kl}}n\\
&=\sum_{ij}T_{ij}(u)\omega^{-ki}\delta_{j,
-l+i}\frac{\omega^{kl}}n=\sum_i \frac{\omega^{k(l-i)}}nT_{i,
-l+i}(u).
\end{align*}

\section{Standard principal representation}

Let $V$ be the fundamental representation of $\frak g\frak l(n)$
with the standard basis $\{ v_i\}_{i=1}^n$. We consider the
endomorphism ring ${\rm End}(V)$ and $E_{ij}\in {\rm End}(V)$:
$$E_{ij}v_k=\delta_{jk}v_i.$$
Using the additive group structure of the index set ${\mathbb
Z}/{n{\mathbb Z}}$, we apply the Fourier transform of the basis
element $v_i$:
$$\varphi_i=\frac{1}{\sqrt n}\sum_{k=1}^n \omega^{ik}v_k.\eqno
(4.1)$$ Under the natural inner product $(v_i|v_j)=\delta_{ij}$ the
vectors $\varphi_i$ satisfy
$$(\varphi_i|\varphi_j)=\delta_{i+j,0}.\eqno (4.2)$$
The transition matrix from the basis $\{\varphi_i\}$ to the basis
$\{ v_i\}$ is given by the inverse Fourier transform:
$$v_i=\frac{1}{\sqrt n}\sum_{k=1}^n \omega^{-ik}\varphi_k.\eqno
(4.3)$$

The action of the principal basis is simply
$$
A_{kl}\varphi_i=\omega^{il}\varphi_{i+k}.
$$

As in the finite dimensional case, the principal realization of the
general linear Lie algebra also provides a representation for the
Yangian. The following result is a reformulation of the evaluation
homomorphism from the Yangian $Y(\frak g\frak l (n))$ to the
universal enveloping algebra $U(\frak g\frak l (n))$.

\begin{proposition} The map $\phi(S_{ij}^{(1)})=A_{ij}$
and $\phi(S_{ij}^{(n)})=0$ for $n\geq 2$ gives rise to a
representation for the Yangian.
\end{proposition}

It is clear that any finite dimensional representation of $U(\frak
g\frak l (n))$ can also be viewed as a representation of the Yangian
as stated in the above proposition.

\section{The representations of the Yangian $Y(sl(3))$ and entangled states}

In this section, we will show that the principal basis given in the
previous sections plays an essential role in the representation
theory of the Yangian $Y(sl(3))$ which has a close relation with the
study of entangled states \cite{KOC} in quantum information.

Let $\lambda_1$ and $\lambda_2$ be the fundamental weights of the
simple Lie algebra $sl(3)$. The irreducible representation
$V(\lambda_2)$ can be viewed as the dual of the irreducible
representation $V(\lambda_1)$.  Suppose $|i\rangle_{1}=u_0, u_1$ and
$u_2$ is the standard basis in $V(\lambda_1)$ (quark states), and
set $|j\rangle_{2}=u_0^*, u_1^*$ and $u_2^*$ are dual base
(antiquark states). Note that $V(\lambda_1)$ is isomorphic to the
3-dimensional vector representation. Let
$|i,j\rangle=|i\rangle_1\otimes |j\rangle_{2}$ ($i,j=0,1,2$) be the
orthonormal basis in the tensor representation $V(\lambda_1)\otimes
V(\lambda_2)$. The $sl(3)$ entangled states with the maximal degree
of entanglement were given in \cite{KOC} as follows.
\begin{equation*}
\psi_j^{(i)}=
\frac{1}{\sqrt{3}}(|0,i-1\rangle+\omega^{i-1}|1,i\rangle+\omega^{2(i-1)}|2,i+2\rangle),
\qquad i,j=1, 2, 3. \eqno(5.1)
\end{equation*}

The Cartan-Weyl basis for the algebra $Y(sl(3))$ are generated by
two sets of generators $I_{\alpha}$ and $J_{\alpha}$. In general any
finite dimensional irreducible representation $V$ of $sl(3)$ can be
lifted to a representation of the Yangian $Y(sl(3))$. In the case of
fundamental representations $V(\lambda_i)$ it is given by the
homomorphism $J_{\alpha}=aI_{\alpha}$ for a fixed constant $a$, and
we denote the resulted presentation by $V(\lambda_i, a)$ . When we
pass the action of the Yangian $Y(sl(3))$ to the tensor product we
need to apply its co-products. According to Drinfeld \cite{D2} the
finite dimensional irreducible representations are classified by
their Drinfeld polynomials, which are determined with help of
Drinfeld's new basis; see also \cite{CP} and \cite{M} for detailed
exposition. We would like to use our principal basis for the tensor
product. First of all the Cartan generators are
$$H_1=E_{11}-E_{22}, H_2=E_{22}-E_{33}, E_{ij}, i\ne j\eqno (5.2)$$
to give the action of the Yangian $Y(sl(3))$ on the tensor product
$V(\lambda_1, a)\otimes V(\lambda_2, b)$. We will give the action of
the Yangian operators $J(H_1), J(H_2), J(E_{ij})$ on the standard
basis of the representations of the Lie algebra $sl(3)$ (given by
Clebsch-Gordon coefficients) according to the following
decomposition of the representations of $sl(3)$
$$
V(\lambda_1)\otimes V(\lambda_2) \cong V(\lambda_1+\lambda_2)\oplus
V(0).\eqno (5.3)$$ If we use the standard basis $|i,j\rangle$ the
action of the Yangian is very complicated and it is not obvious to
see whether there exist certain rules.

We use the entangled states $\psi_j^{(m)}$ given by Eqs. (5.1)-(5.3)
as the basis for the representation space, and slightly modify the
principal basis $A_{ij}$ ($ij\ne00$) as follows.

$$T_i^{(j)}=\omega^{3-i+1}A_{i-1,j-1},\;\; i,j=1,2,3\;\;{\rm and}\;\; {\rm for}\;\;j=1, i\ne
1.\eqno (5.4)$$ Then comparing with the basis given by Eq. (5.2), we
know that
$$T_2^{(1)}=H_1-\omega^2 H_2,\;\; T_3^{(1)}=H_1-\omega H_2,$$
$$T_1^{(2)}=E_{12}+E_{23}+E_{31},\;T_2^{(2)}=E_{12}+\omega E_{23}+\omega^2 E_{31},\;
T_3^{(2)}=E_{12}+\omega^2 E_{23}+\omega E_{31},$$
$$T_1^{(3)}=E_{13}+E_{21}+E_{32},\;T_2^{(3)}=E_{13}+\omega E_{21}+\omega^2 E_{32},\;
T_3^{(3)}=E_{13}+\omega^2 E_{21}+\omega E_{32}.\eqno (5.5)$$ In
terms of Eq. (5.5) the Yangian operators $J(T_2^{(1)}),
J(T_3^{(1)}), J(T_1^{(2)}),\cdots$ can be defined according to the
action of Yangian operators expressed in terms of $H_1, H_2$ and
$E_{ij}(i\ne j)$.

A direct computation proves the following result.

{\bf Theorem 5.1}\quad The action of the Yangian $Y(sl(3))$ on the
tensor product of $V(\lambda_1, a)\otimes V(\lambda_2, b)$ can be
expressed in the action of the principal basis $J(T_i^{(j)})$ on the
entangled states $\psi_k^{(m)} (k,m=1,2,3)$ with the maximal degree
of entanglement in a simple way. The explicit action is given by the
following equation:
$$J(T_i^{(j)})\psi_k^{(m)}=[a\omega^{(j-1)(k-1)}-b\omega^{(i-1)(m-1)}
+\frac{3}{2}\delta_{i+k-1,1}\delta_{m+j-1,1}\omega^{(j-1)(k-1)}
-\frac{3}{2}\delta_{k,1}\delta_{m,1}]\psi_{i+k-1}^{(m+j-1)}.$$

We remark that the vectors $\psi_k^{(m)}$ provide a nice basis which
simplifies the action of $Y(sl(3))$ dramatically, and the transition
operators for the entangled states $\psi_k^{(m)} (k,m=1,2,3)$ are
with the maximal degree of entanglement in quantum computation.

From Theorem 5.1 it is easy to get the following conclusion
(\cite{CP}, \cite{M}).

{\bf Corollary 5.2}\quad Let $W=V(\lambda_1, a)\otimes V(\lambda_2,
b)$. If $|a-b|\ne \frac{3}{2}$, then $W$ is an irreducible
 representation of the Yangian $Y(sl(3))$. Otherwise, $W$
has a unique proper $Y(sl(3))$-subrepresentation $V$ given as
follows.

(1) If $a-b=\frac{3}{2}$, we have $V\cong V(0,0)$, and $W/V\cong
V(\lambda_1+\lambda_2)$ as vector spaces.

(2) If $a-b=-\frac{3}{2}$, we have $V\cong V(\lambda_1+\lambda_2)$
as vector spaces.

\section*{Acknowledgments} Jing is grateful to
the support of NSA grant H98230-06-1-0083 and NSFC's Overseas
Distinguished Youth Grant (10728102). This work was supported in
part by the National Natural Science Foundation of China (10575053,
10571091, 10621101), NKBRPC (2006CB805905), Program for New Century
Excellent Talents in University.

\bibliographystyle{amsalpha}

\end{document}